\documentclass{amsart}
\usepackage{amsfonts}
\usepackage{amsmath}
\usepackage{amssymb}

\theoremstyle{plain}
\newtheorem{theorem}{Theorem}[section]

\newtheorem{proposition}[theorem]{Proposition}
\theoremstyle{definition}
\newtheorem{definition}[theorem]{Definition}

\theoremstyle{remark}

\input{tcilatex}

\begin{document}
\title{The Courant type algebroids, the coadjoint orbits and related integrable
ows}
\author{Anatolij K. Prykarpatski$^{1,2,3)}$ \ and Victor A. Bovdi$^{2)}$}
\address{$^{1)}$ Cracow University of Technology, Cracow, Poland}
\email{pryk.anat@cybergal.com}
\address{ $^{2)}$ Department of Mathematics at the UAEU, \ Al Ain, UAE}
\email{vbovdi@gmail.com}
\address{ $^{3)}$ The Lviv Polytechnic National University, Lviv, 79000,
Ukraine}
\thanks{Authors declare: No conflict of interest}
\maketitle

\begin{abstract}
Poisson structures related with the affine Courant type algebroid are
analyzed, including \ those related with cotangent bundles on Lie group
manifolds. A special attantion is paid to Courant type algebroids and
related R-structures \ on them, generated by suitably defined tensor
mappings. \ There are constructed Lie-Poisson brackets invariant with
respect to the coadjoint action of the loop diffeomorphisms group and
described the related Courant type algebroids. \ The corresponding
integrable Hamiltonian flows, generated by Casimir functionals and
generalizing the so called heavenly type differential systems, describing
diverse geometric structures of conformal type on finite dimensional
Rieamnnian manifolds are described.
\end{abstract}

\section{Introduction}

The Lie algebroid \ \ \cite{Mack} \ as a mathematical object is an \textit{%
unrecognized} part of the folklore of differential geometry. They have been
introduced repeatedly into differential geometry since the early 1950's, and
also into physics and algebra, under a wide variety of names, chiefly as
infinitesimal invariants associated to geometric structures: in connection
theory, as a means of treating de Rham cohomology by algebraic methods, as
invariants of foliations and pseudogroups of various types, in symplectic
and Poisson geometry and, in a more algebraic setting, as algebras of
differential operators associated with vector bundles and with infinitesimal
actions of Lie groups. The algebroid structures found recently diverse
applications in geometry of Poisson \ \cite{GaVaVo} and Lagrangian manifolds
\ \ \cite{LeMaMa}, in mechanical sciences \ \ \cite{JiLe} and other branches
of modern applied research.\

\begin{definition}
Let $M$ be a manifold. A Lie algebroid $(E;[[.,.]],\rho ,M)$ on $M,$ or with
base $M,$ is a vector bundle $E$ $\rightarrow $ $M,$ together with a bracket
$[[.,.]]:\Gamma (E)\times \Gamma (E)\rightarrow \Gamma (E)$ on the module $%
\Gamma (E)$ of global sections of $E$, and a vector bundle morphism $\rho :$
$E\rightarrow T(M)\ $ from $E$ to the tangent bundle $T(M)$ of $M,$ called
the anchor of $E,$ such that
\end{definition}

(i) the bracket on $\Gamma (E)$ is $\mathbb{R}$-bilinear, skew-symmetric and
satisfies the Jacobi identity;

(ii) $[[\alpha ,f\beta ]]=f[[\alpha ,\beta ]]+\ \rho (\alpha )f$ $\beta $
for all $\alpha ,\beta \in \Gamma (E)$ and all smooth functions $\ \ f\in
\mathcal{D}(M);$

(iii) $\rho ([[\alpha ,\beta ]])=[\rho (\alpha ),\rho (\beta )]$ \ for all $%
\ \ \alpha ,\beta \in \Gamma (E).$

The anchor $\rho $ ties here the bracket on $\Gamma (E)$ to the vector field
structure on $M$ as a module over $\mathcal{D}(M),$ the algebra of smooth
functions $f:M\rightarrow \mathbb{R}.$

Consider now the product $T(M)\ltimes T^{\ast }(M)$ of the tangent $T(M)$
and its cotangent $T^{\ast }(M)$ bundles over the manifold $M.$ Then the
\textit{canonical Courant bracket} \ \ \cite{Cour} \ on the $\mathcal{D}(M)-$
module $\mathcal{A}(M):=T^{\ast }(M)\times T(M)\simeq \left( T(M)\times
T^{\ast }(M)\right) ^{\ast }$ is defined as
\begin{equation}
\lbrack \lbrack (\alpha ,a),(\beta ,b)]]:=(L_{a}\beta -L_{b}\alpha +\frac{1}{%
2}d(\alpha (b)-\beta (a)),\ [a,b])  \label{A1}
\end{equation}%
for any $(\alpha ,a),(\beta ,b)\in \mathcal{\ }T^{\ast }(M)\times T(M),$
satisfying \ \ \cite{AbMa,Godb,KoMiSl} \ the usual Jacobi identity.

\begin{definition}
The bundle $\mathcal{A}(M)=T^{\ast }(M)\times T(M)\ $ \ jointly with the
bracket \ (\ref{A1}) and the natural morphism projection mapping $\rho :%
\mathcal{A}(M)\rightarrow T(M)$ is called the Courant algebroid.
\end{definition}

Let now assume that the cotangent space \ $T^{\ast }(M)$ is endowed with its
own Poisson structure $P:T^{\ast }(M)\rightarrow T(M).$ \ Then having put,
by definition, $a:=P\alpha ,b:=P\beta \in \mathcal{T}(M),$ one can easily
observe that the Courant bracket \ (\ref{A1}) becomes in its second term
identity, reducing to the next \ bracket on the cotangent space $T^{\ast
}(M):$
\begin{eqnarray}
\lbrack \lbrack \alpha ,\beta ]] &=&L_{P\alpha }\beta -L_{P\beta }\alpha +%
\frac{1}{2}d(\alpha (P\beta )-\beta (P\alpha )=  \notag \\
&&  \notag \\
&=&i_{P\beta }d\alpha -i_{P\alpha }d\beta -\frac{1}{2}d(\alpha (P\beta
)-\beta (P\alpha )  \label{A4}
\end{eqnarray}%
for any $\alpha ,\beta \in \Lambda ^{1}(M),$ \ satisfying \ the Jacobi
identity. Thus, the triple $(T^{\ast }(M);[[\cdot ,\cdot ]],P)$ becomes a
Lie algebroid with the \textit{ancor} $\ P:T^{\ast }(M)\rightarrow T(M),$
considered as a Lie algebra morphism
\begin{equation}
P[[\alpha ,\beta ]]=[P\alpha ,P\beta ],  \label{A5}
\end{equation}%
satisfied for any $\alpha ,\beta \in T^{\ast }(M)\simeq \Lambda ^{1}(M).$ \

As an instructive example of the construction above we consider a semisimple
Lie group $G$ and its Lie algebra $\mathcal{G}\simeq T_{e}(G)$ at the unity
element $e\in G,$ consisting of the left invariant vector fields on $G.$
Assume that \ \ \cite{AbMa,Arno,KoMiSl,BlPrSa} $\Omega :\mathcal{G}%
\rightarrow \mathcal{G}^{\ast }$ is a symplectic structure on $G,$ allowing
to construct for any $\alpha ,\beta \in \mathcal{G}^{\ast }$ the adjoint
left-invariant vector fields as $X_{\alpha }:=\Omega ^{-1}\alpha ,$ $%
X_{\beta }:=\Omega ^{-1}\beta \in \mathcal{G},$ subject to which the related
Poisson bracket
\begin{equation}
\lbrack \lbrack \alpha ,\beta ]]:=i_{[X_{\alpha },X_{\beta }]}\Omega
\label{A5a}
\end{equation}%
satisfies the Jacobi identity. The latter, in particular, means that the
constructed object $(\mathcal{G}^{\ast };[[.,.]],\Omega ^{-1},G)$ is also a
reduced Lie algebroid. Moreover, the Lie bracket \ (\ref{A5a}), owing to the
Cartan representation of the Lie derivative $L_{X}=i_{X}d+di_{X},$ $X\in
\mathcal{G},$ \ \ \cite{AbMa,Arno,Godb,KoMiSl}, \ can be rewritten as
\begin{equation}
\begin{array}{c}
\lbrack \lbrack \alpha ,\beta ]](Z)=\left( i_{[X_{\alpha },X_{\beta
}]}\Omega \right) (Z)=[L_{X_{\alpha }},i_{X_{\beta }}]\Omega (Z)= \\
=L_{X_{\alpha }}i_{X_{\beta }}\Omega (Z)-i_{X_{\beta }}L_{X_{\alpha }}\Omega
(Z)= \\
=\ i_{X_{\alpha }}di_{X_{\beta }}\Omega (Z)+d\left( i_{X_{\alpha
}}i_{X_{\beta }}\Omega \right) (Z)-i_{X_{\beta }}i_{X_{\alpha }}d\Omega
(Z)-i_{X_{\beta }}di_{X_{\alpha }}\Omega (Z)= \\
=\ i_{X_{\alpha }}di_{X_{\beta }}\Omega (Z)-i_{X_{\beta }}di_{X_{\alpha
}}\Omega (Z)+d\left( \Omega (X_{\beta },X_{\alpha })\right) (Z)= \\
=X_{\alpha }\Omega (X_{\beta },Z)-Z\Omega (X_{\beta },X_{\alpha })-\Omega
(X_{\beta },[X_{\alpha },Z])- \\
-X_{\beta }\Omega (X_{\alpha },Z)+Z\Omega (X_{\alpha },X_{\beta })+ \\
+\Omega (X_{\alpha },[X_{\beta },Z])+d\left( \Omega (X_{\beta },X_{\alpha
})\right) (Z)= \\
=-\Omega (X_{\beta },[X_{\alpha },Z])+\Omega (X_{\alpha },[X_{\beta
},Z])-d\left( \Omega (X_{\alpha },X_{\beta })\right) (Z) \\
=\left( ad_{X_{\beta }}^{\ast }\left( i_{X_{\alpha }}\Omega \right)
-ad_{X_{\alpha }}^{\ast }\left( i_{X_{\beta }}\Omega \right) \right)
(Z)-d\left( \Omega (X_{\alpha },X_{\beta })\right) (Z)= \\
=\left( ad_{\Omega ^{-1}\beta }^{\ast }\alpha -ad_{\Omega ^{-1}\alpha
}^{\ast }\beta \right) (Z)-1/2d\left( \alpha (\Omega ^{-1}\beta )-\beta
(\Omega ^{-1}\alpha )\right) (Z),%
\end{array}
\label{A5c}
\end{equation}%
where we made use of the invariance conditions $Z\alpha (X)=0=Z\beta (X)$
for arbitrary $\alpha ,\beta \in \mathcal{G}^{\ast }$ and $X,Z\in \mathcal{G}%
,$ as well as we denoted by $ad^{\ast }:\mathcal{G\times G}^{\ast
}\rightarrow \mathcal{G}^{\ast }$ the natural coadjoint action of the Lie
algebra $\mathcal{G}$ on the adjoint space $\mathcal{G}^{\ast }.$ \ \ The \
obtained expression \ (\ref{A5c}) on $\mathcal{G}^{\ast }$ can be rewritten
as
\begin{equation}
\lbrack \lbrack \alpha ,\beta ]]=ad_{\rho (\beta )}^{\ast }\alpha -ad_{\rho
(\alpha )}^{\ast }\beta -1/2\ d\left[ \alpha (\rho (\beta ))-\beta (\rho
(\alpha ))\right] ),  \label{A5d}
\end{equation}%
where $\rho =\Omega ^{-1}:\mathcal{G}^{\ast }\rightarrow \mathcal{G}$
denotes the corresponding anchor mapping subject to the reduced Courant Lie
algebroid $(\mathcal{G}^{\ast };[[.,.]],\rho ,G).$ Assume now that we are
given a Lie-algebroid $(\mathcal{G}_{h}^{\ast };[[.,.]],\rho _{h},G_{h}),$
whose anchor $\rho _{h}:\mathcal{G}^{\ast }\rightarrow \mathcal{G}$ is a Lie
algebra homomorfism, being not necesary related to a symplectic structure $\
$on $G_{h}$ and \textit{a priori} satisfying the Poisson bracket \ (\ref{A5d}%
). Then our inverse problem consists in describing at least sufficient
conditions on the anchor $\rho _{h}:\mathcal{G}^{\ast }\rightarrow \mathcal{G%
},$ under which the bracket \ (\ref{A5d}) will satisfy the Jacobi condition.

As a simple guiding construction for our proceeding Courant algebroid
constructions, let us consider the cohomology group $H^{1}(G;\mathbb{C})$ of
the Lie group $G$ and observe that for any $\ \{\alpha \},\ \{\beta \}\in
H^{1}(G;\mathbb{C}),\ $\ $\alpha ,\beta \in \mathcal{G}^{\ast },$ the
Poisson bracket \ \ (\ref{A5d}) reduces to the Poisson bracket
\begin{equation}
\begin{array}{c}
\lbrack \lbrack \{\alpha \},\{\beta \}]]=\ \{ad_{\rho (\beta )}^{\ast
}\alpha -ad_{\rho (\alpha )}^{\ast }\beta \},%
\end{array}
\label{A5e}
\end{equation}%
on $H^{1}(G;\mathbb{C}),$ satisfying the Jacobi condition. Thus, the
cohomology group $H^{1}(G;\mathbb{C})$ is simultaneously the Lie algebra
with respect \ to the Lie product \ (\ref{A5e}), satisfying the induced
property $\Omega ^{-1}\rightarrow \Omega _{h}^{-1}:H^{1}(G;\mathbb{C}%
)\rightarrow $ $\mathcal{G}/\mathcal{H}\simeq T(G/H),$ where \ $H\subset G$
\ denotes the normal \textit{Hamiltonian subgroup} of $G,$ whose Lie algebra
$\mathcal{H\subset G}$ consists of the Hamiltonian shifts $\Omega (h)\in
\mathcal{H}\ $for all closed elements $h\in \mathcal{G}^{\ast },$ $dh=0.$\
The latter makes it possible to construct the reduced Lie algebroid: $%
(T^{\ast }(G_{h});[[.,.]]_{h},\rho _{h},G_{h})\ $ with the Lie bracket
\begin{equation}
\lbrack \lbrack \tilde{\alpha},\tilde{\beta}]]_{h}=ad_{\rho _{h}(\tilde{\beta%
})}^{\ast }\tilde{\alpha}-ad_{\rho _{h}\left( \tilde{\alpha}\right) }^{\ast }%
\tilde{\beta}  \label{A5f}
\end{equation}%
for any $\tilde{\alpha},\tilde{\beta}\in T^{\ast }(G_{h}),$ where $%
G_{h}:=G/H\ $ and the anchor $\rho _{h}:=\Omega _{h}^{-1}:T^{\ast
}(G_{h})\rightarrow $ $\mathcal{G}/\mathcal{H}.$

\section{Courant type algebroids and the related R-structures}

In what follows futrther we \ will deal with the semi-direct product bundle
Lie algebra $\mathcal{A}^{\ast }(M):=T(M)\ltimes T^{\ast }(M),$ whose Lie
product is defined as
\begin{equation}
\lbrack (a,\alpha ),(b,\beta )]:=([a,b],ad_{b}^{\ast }\alpha -ad_{a}^{\ast
}\beta )  \label{A1a}
\end{equation}%
for any $(a,\alpha ),(b,\beta )\in \ \mathcal{A}^{\ast }(M),$ satisfying the
Jacobi identity. Moreover, the Lie algebra \ $\mathcal{A}^{\ast }(M)$ proves
to be metrized \-\ \ \cite{FaTa,Seme} \ with respect to the dual $ad$%
-invariant pairing:%
\begin{equation}
((a,\alpha )|(b,\beta )):=\alpha (b)+\beta (a),  \label{A2}
\end{equation}%
as for arbitrary $(\alpha ,a),(\beta ,b)$ and $(c,\gamma )\in \mathcal{A}%
^{\ast }(M)$ the identity%
\begin{equation}
\left( \lbrack (a,\alpha ),(b,\beta )]|(c,\gamma )\right) =(a,\alpha
)|[(b,\beta ),(c,\gamma )])  \label{A3}
\end{equation}%
holds. Take now into account that the cotangent space $\mathcal{A}(M)^{\ast
}\simeq T(M)\times T^{\ast }(M)$ posessess the canonical Lie-Poisson
structure, defined by means of the following \ bracket
\begin{equation}
\lbrack \lbrack ((l,p)|X),((l,p)|Y)]]_{Lie}=((l,p)|[X,Y]),  \label{A2a}
\end{equation}%
for all $X,Y\in \mathcal{A}(M)^{\ast }$ and any fixed element $(l,p)\in
\mathcal{A}(M).$ \

To construct a Courant type algebroid $\mathcal{A}(M)$ let us take a tensor
element $r\in \mathcal{A}^{\ast }(M)\otimes \mathcal{A}^{\ast }(M)$ jointly
with the related linear mapping $\ $ $r:\mathcal{A}(M)\rightarrow \mathcal{A}%
^{\ast }(M)\ $ and define on the bundle $\mathcal{A}(M)$ the following
bracket
\begin{equation}
\lbrack \lbrack (\alpha ,a),(\beta ,b)]]_{r}:=ad_{r((\beta ,b))}^{\ast }\
(\alpha ,a)-ad_{r((\alpha ,a))}^{\ast }\ (\beta ,b)  \label{A3a}
\end{equation}%
for any $(\alpha ,a),(\beta ,b)\in $ $\mathcal{A}(M).$ \ The following
proposition, describing the conditions to be imposed on the tensor element $%
r\in \mathcal{A}^{\ast }(M)\otimes \mathcal{A}^{\ast }(M),$ \ holds.

\begin{proposition}
Let a tensor element $r\in \mathcal{A}^{\ast }(M)\otimes \mathcal{A}^{\ast
}(M)$ allow the direct sum splitting : $r=k\oplus \eta ,$ where $k\in
\mathcal{A}^{\ast }(M)\overset{Sym}{\otimes }\mathcal{A}^{\ast }(M)$ \ is
its symmetric and $\eta \in \mathcal{A}^{\ast }(M)\wedge \mathcal{A}^{\ast
}(M)$ antisymmetric non-degenerate parts. If, in addition, \ the related
mapping $D:=\ k\circ \eta ^{-1}:\mathcal{A}^{\ast }(M)\rightarrow \mathcal{A}%
^{\ast }(M)$ is a derivation of the algebra $\mathcal{A}^{\ast }(M),$ that
is
\begin{equation}
D[(a,\alpha ),(b,\beta )]=[D(a,\alpha ),(b,\beta )]+[(a,\alpha ),D(b,\beta )]
\label{A9}
\end{equation}%
for any $(a,\alpha ),(b,\beta )\in $ $\mathcal{A}^{\ast }(M),$ then the
mapping $\rho _{r}=pr_{1}\circ r:$ $\mathcal{A}(M)\rightarrow T(M)$ \
\begin{equation}
\rho _{r}\ [[(\alpha ,a),(\beta ,b)]]_{r}=[\rho _{r}(\alpha ,a),\rho
_{r}(\beta ,b)]  \label{A10}
\end{equation}%
defines a Lie algebra homomorphism and the triple $(\mathcal{A}%
(M);[[.,.]],\rho _{r})$ is a Courant type algebroid.
\end{proposition}

In addition, the following important corollaries are infered from reasonings
above: \textbf{\ }if the invertible tensor $r\in \mathcal{A}(M)\otimes
\mathcal{A}(M)\ \ $satisfies the conditions, formulated above, \ then the
corresponding triple $(\mathcal{A}(M)^{\ast },[[\cdot ,\cdot ]],\rho )$\ is
a generalized Courant type algebroid with the related ancor morphism $\rho :%
\mathcal{A}(M)^{\ast }\rightarrow T(M),$ defined by means of the composition
$\rho =pr_{1}\circ r$ \ with the projection mapping $pr_{1}:T(M)\times
T^{\ast }(M)\rightarrow T(M)$ on the first component. Moreover, the
following deformed commutator structure
\begin{equation}
\lbrack \lbrack (a,\alpha ),(b,\beta )]]_{R}:=[[R(a,\alpha ),(b,\beta
)]]+[[(a,\alpha ),R(b,\beta )]]  \label{A11}
\end{equation}%
determines for any $(a,\alpha ),(b,\beta )\in \ \mathcal{A}(M)$\ the
corresponding $R-$structure on $\mathcal{A}(M),$ $R=D^{-1},$ subject to
which the second Lie bracket $[[\cdot ,\cdot ]]_{D}$ on $\mathcal{A}(M)\ $\
also satisfies the Jacobi identity. The latter makes it possible to
construct on the space $\mathcal{A}(M)^{\ast }$ the deformed Lie-Poisson
bracket
\begin{equation}
\lbrack \lbrack ((l,p)|X),((l,p)|Y))]]_{R}=((l,p)|[[X,Y]]_{R})\ \
\label{A12}
\end{equation}%
for all $X,Y\in \mathcal{A}(M)$ and any fixed element $(l,p)\in \mathcal{A}%
(M)^{\ast },$ subject to which the Casimir functionals of the Lie bracket \ (%
\ref{A2a}) naturally generate \ \ \cite{BlPrSa,Blas-1,Blas-2,FaTa,Seme} \ an
infinite hierarchy of commuting to each other completely integrable
Hamiltonian systems.

\section{Remarks on the Courant type algebroid foliation and related
Hamiltonian flows}

Let $(E;[[\cdot ,\cdot ]],\rho ),E:=\mathcal{A}(M),$ \ be a Courant type
algebroid over a manifold $M,$ for which the characteristic space $\rho
(E_{x})\subset T_{x}(M),x\in M,$ is involutive and finitely-generated. \ On
the algebroid $(E;[[\cdot ,\cdot ]],\rho )$ there is naturally defined the
external differentialof $d_{E}:\Gamma (\Lambda (E^{\ast }))\rightarrow
\Gamma (\Lambda (E^{\ast })),$ where $\Lambda (E^{\ast }):=\oplus _{k\in
\mathbb{Z}_{+}}$ $\wedge ^{k}(E^{\ast }),$ as follows:%
\begin{equation}
\begin{array}{c}
\left( d_{E}\alpha ^{(k)}\right)
(A_{0},A_{1},...,A_{k}):=\sum_{i=0}^{k}(-1)^{i}\rho (A_{i})_{E}\alpha
^{(k)}((A_{0},A_{1},...,\hat{A}_{i},...,A_{k})+ \\
\\
+\sum_{i<j=0}^{k}(-1)^{i+j}\alpha ^{(k)}([[A_{i},A_{j}]],A_{0},A_{1},...,%
\hat{A}_{i},...,\hat{A}_{j},...,A_{k}),%
\end{array}
\label{B1}
\end{equation}%
for $\alpha ^{(k)}\in \Gamma (\Lambda ^{k}(E^{\ast }))\ $and arbitrary $%
A_{i}\in \Gamma (E),i=\overline{0,k},$ and $k\in \mathbb{Z}_{+},$ satisfying
the natural algebraic cohomology complex property $d_{E}d_{E}=0.$ The
differential \ (\ref{B1}) makes it possible to determine also the Lie
derivative
\begin{equation}
\mathcal{L}_{A}^{(E)}:=i_{A}d_{E}+d_{E}i_{A}  \label{B2}
\end{equation}%
along a section $A\in \Gamma (E).$ It is obviously that for $E=T^{\ast
}(M)\times T(M),\rho :=$ $pr_{T(M)}$ and $A\in \Gamma (T^{\ast }(M)\times
T(M)),$ the external differential $d_{E}$ transits into the usual
differential $\ d\circ pr_{\Lambda (M)}:$ $T(M)\times \Lambda (M)\rightarrow
\Lambda (M).$

Let a two-form $\omega ^{(2)}\in \Gamma (\Lambda ^{2}(E))$ be closed and
invariant with respect to a vector $K\in \Gamma (E),$ that is
\begin{subequations}
\begin{equation}
d_{E}\omega ^{(2)}=0,\text{ \ \ }\mathcal{L}_{K}^{(E)}\omega ^{(2)}=0
\label{B3}
\end{equation}%
on $M.$ The latter\ \ \ simply means that there exists a locally defined
smooth function $H_{1}:\Gamma (E^{\ast })\rightarrow \mathbb{R},$ such that
\end{subequations}
\begin{equation}
i_{K}\omega ^{(2)}=-d_{E}H_{1}.  \label{B4}
\end{equation}%
If the function $H_{1}:\ \Gamma (E^{\ast })\rightarrow \mathbb{R}$ is
defined globally, then the vector $K\in \Gamma (E)$ will be called a \textit{%
Hamiltonian flow} on the bundle $E.$

Consider now a coordinate vector set $e:=\{e^{_{i}}\in \Gamma (E^{\ast }):i=%
\overline{1,m}\}$ on the vector bandle $\Gamma (E),$ the corresponding basis
of its differenrials $\{d_{E}e^{i}\in \Gamma (T^{\ast }(E^{\ast })):i=%
\overline{1,m}\}\ \ $ and take a linear and invertible subject to the second
component mapping $Q_{E}:\Gamma (E^{\ast })\times \Gamma (T^{\ast }(E^{\ast
}))\rightarrow \Gamma (T^{\ast }(E^{\ast })),$ determining the elements
\begin{equation}
d_{E}^{\ast }e=\ \sum_{j=\overline{1,m}}Q_{E}(e)_{j}d_{E}e^{j}\in \Gamma
(T^{\ast }(E^{\ast }))\   \label{B5}
\end{equation}%
by means of some mappings $Q_{E}(e)_{j}\in \mathrm{End}$ $\Gamma (E^{\ast
}), $ $j=\overline{1,m}.$ \ The expression \ (\ref{B5}) makes \ it \ \ \
possible to define a second external differential $d_{E}^{\ast }:\Gamma
(\Lambda (E^{\ast }))\rightarrow \Gamma (\Lambda (E^{\ast })),$ satisfying
the algebraic cohomology complex property $\ d_{E}^{\ast }d_{E}^{\ast }=0\ $%
\ and anticommuting to the external deformed differential $d_{E}:\Gamma
(\Lambda (E^{\ast }))\rightarrow \Gamma (\Lambda (E^{\ast })):$
\begin{equation}
d_{E}d_{E}^{\ast }+d_{E}^{\ast }d_{E}=0.  \label{B6}
\end{equation}%
As the symplectic two-form $\omega ^{(2)}\in \Gamma (\Lambda ^{2}(E^{\ast
})) $ is, by definition, $d_{E}$ - closed, $d_{E}\omega ^{(2)}=0,$ \ \ we
will for furthere assume additionally also its $d_{E}^{\ast }$-closedness,
that is $d_{E}^{\ast }\omega ^{(2)}=0.$ \

\begin{definition}
External derivations $d_{E},d_{E}^{\ast }:\Gamma (\Lambda (E^{\ast
}))\rightarrow \Gamma (\Lambda (E^{\ast })),$ \ satisfying the properties $%
d_{E}d_{E}=0,$ $\ d_{E}^{\ast }d_{E}^{\ast }=0$ and $d_{E}d_{E}^{\ast
}+d_{E}^{\ast }d_{E}=0,$ are called compatible.
\end{definition}

Impose additionally on the symplectic two-form $\omega ^{(2)}\in \Gamma
(\Lambda ^{2}(E^{\ast }))$ its $K$ - ivariance:
\begin{equation}
\ \mathcal{L}_{K}^{(E,\ast )}\omega ^{(2)}:=0\ \   \label{B2aa}
\end{equation}%
with respect to the related Lie-derivative
\begin{equation}
\mathcal{L}_{K}^{(E,\ast )}:=i_{K}d_{E}^{\ast }+d_{E}^{\ast }i_{K}.
\label{B2a}
\end{equation}%
Whence one easily derives that the following expressions
\begin{equation}
\ i_{K}\omega ^{(2)}=-d_{E}^{\ast }H_{1}=-d_{E}H_{2}  \label{B6a}
\end{equation}%
holds for some smooth mapping $H_{1},H_{2}:\Gamma (E^{\ast })\rightarrow
\mathbb{R}.$ $\ $Moreover, multiplying the right hand side of\ the equality $%
\ -d_{E}^{\ast }H_{1}=-d_{E}H_{2}$ by the differential $\ d_{E}^{\ast
}:\Gamma (\Lambda (E^{\ast }))\rightarrow \Gamma (\Lambda (E^{\ast })),$ one
easily obtains that
\begin{equation}
\ \ -d_{E}^{\ast }d_{E}^{\ast }H_{1}=d_{E}(d_{E}^{\ast }H_{2})=0,
\label{B6b}
\end{equation}%
meaning that there exists a smooth mapping $H_{3}:\Gamma (E^{\ast
})\rightarrow \mathbb{R},$ such that
\begin{subequations}
\begin{equation}
d_{E}^{\ast }H_{2}=d_{E}H_{3}.  \label{B6c}
\end{equation}%
The obtained above relationship \ (\ref{B6a}) can be recurrently continued
and equivalently rewritten as the modified \textit{Lenard type} mapping $%
Q_{E}^{(\ast )}:\Gamma (T^{\ast }(E^{\ast }))\rightarrow \Gamma (T^{\ast
}(E^{\ast })),$ acting as
\end{subequations}
\begin{equation}
d_{E}H_{j+1}:=Q_{E}^{(\ast )}\left( d_{E}H_{j}\right) =d_{E}^{\ast }H_{j}
\label{B6e}
\end{equation}%
for the \ set of smooth mappings $H_{j}:\Gamma (E^{\ast })\rightarrow
\mathbb{R},$ $j=\overline{1,m},$ $H_{2}=H,d_{E}^{\ast }H_{m}=0,$ generated
by the $\ $second external differential $d_{E}^{\ast }:\Gamma (\Lambda
(E^{\ast })\rightarrow \Gamma (\Lambda (E^{\ast }),$ defined by the
expressions \ (\ref{B5}).

\ It is also easy to observe that there exists the one-form $\beta
_{1}^{(1)}\in \Gamma (\Lambda ^{1}(E)),$ satisfying the condition $%
d_{E}\beta _{1}^{(1)}=\omega _{1}^{(2)}:=\omega ^{(2)}\ $and generating \
the second closed 2 - form $\omega _{2}^{(2)}:=d_{E}\circ Q_{E}^{(\ast
)}(\beta _{1}^{(1)})\in $ $\Gamma (\Lambda ^{2}(E^{\ast })),$ which will be
both $d_{E}$ - and $d_{E}^{\ast }$ - closed too. \ The latter makes it
possible to construct a countable hierarchy of symplectic structures $\omega
_{j}^{(2)}\in \Gamma (\Lambda ^{2}(E^{\ast })),j=\overline{1,m},$ such that $%
\omega _{j+1}^{(2)}=d_{E}\circ Q_{E}(\beta _{j}^{(1)}),\omega
_{j}^{(2)}:=d_{E}\beta _{j}^{(2)}$ for $j=\overline{1,m},$ and satisfying
the next recurrent relationship
\begin{equation}
\omega _{j+1}^{(2)}=d_{E}\circ Q_{E}^{(\ast )}\circ d_{E}^{-1}(\omega
_{j}^{(2)}),\text{ \ }d_{E}\circ Q_{E}^{(\ast )}(\beta _{m}^{(1)})=0,
\label{B7}
\end{equation}%
$\ $jointly with the next countable hierarchy of differential relationships
\begin{equation}
i_{K}\omega _{j}^{(2)}=-d_{E}H_{j+1}  \label{B8}
\end{equation}%
for all integers $j=\overline{1,m}.$ Moreover, one easily checks that all
constructed above Hamiltonian functions $H_{j}:\Gamma (E^{\ast })\rightarrow
\mathbb{R},$ $j=\overline{1,m},$ are commuting to each other, that is
\begin{equation}
\{\{H_{j},H_{i}\}\}_{s}=0  \label{B9}
\end{equation}%
for all $i,j=\overline{1,m}\ $ with respect to the following Poisson
brackets:
\begin{equation}
\{\{f,g\}\}_{s}:=\omega _{s}^{(2)}(K_{f},K_{g}),\text{ }i_{K_{f}}\omega
_{s}^{(2)}=-d_{E}f,i_{K_{f}}\omega _{s}^{(2)}=-d_{E}f,  \label{B10}
\end{equation}%
defined for $s=\overline{1,m},$ and arbirary smooth functions $f,g:\Gamma
(E^{\ast })\rightarrow \mathbb{R}.$ The described above properties of the
Courant type algebroid $\ $foliation over the manifold $M$ can be summarized
as the next theorem.

\begin{theorem}
The Courant type algebroid foliation $(E;[[\cdot ,\cdot ]],\rho ),$ $E=%
\mathcal{A}(M),$ \ endowed with two compatible external differentials $%
d_{E},d_{E}^{\ast }:\Gamma (\Lambda (E^{\ast }))\rightarrow \Gamma (\Lambda
(E^{\ast })),\ $ generate \ a finite set of commuting to each other
Hamiltonian flows $K_{j}:\Gamma (E^{\ast })\rightarrow \Gamma (E),j=%
\overline{1,m},$ $\ $ $\ $on$\ $a coordinate set $e$ $\in \Gamma (E^{\ast })
$ of the bundle $\Gamma (E),$ that is \
\begin{equation}
\ [[K_{j},K_{i}]]=0  \label{B11}
\end{equation}
for all $i,j=\overline{1,m},$ where
\begin{equation}
K_{j}(e)=\{\{H_{1},e\}\}_{_{j}}\   \label{B12}
\end{equation}
and $K_{1}=K$ $\in \Gamma (E).$
\end{theorem}

It is here worthy to observe, that the external deformed differential $%
d_{E}^{\ast }:\Gamma (\Lambda (E^{\ast }))\rightarrow \Gamma (\Lambda
(E^{\ast })),\ $defined by means of the relationship \ (\ref{B5}) and
satisfying the constraint $d_{E}^{\ast }d_{E}^{\ast }=0,$ can be interpreted
as a flat connection $d_{E}^{\ast }:\Gamma (E^{\ast })\rightarrow T^{\ast
}(E^{\ast })\otimes \Gamma (E^{\ast })$ on $\Gamma (E), $ acting as%
\begin{eqnarray}
d_{E}^{\ast }A\ &=&\sum_{j=\overline{1,m}}[d_{E}(e_{j}(e))e^{j}(A)+e_{j}(e)%
\left( d_{E}e^{j}\right) (A)]+  \notag \\
&&+\sum_{i,j=\overline{1,m}}\ e_{i}(e)Q_{E}(e)_{j}^{i}\left(
d_{E}e^{j}\right) (A)  \label{B13}
\end{eqnarray}%
on any vector $A:=\sum_{j=1}^{m}e_{j}(e)\ A^{j}\in \Gamma (E),$ where $%
A^{j}:=e^{j}(A),j=\overline{1,m}.\ $Really, as it is easy to check, the
curvature 2-form
\begin{equation}
\Omega ^{(2)}:=d_{E}^{\ast }d_{E}^{\ast }=d_{E}\left( \sum_{j=\overline{1,m}%
}Q_{E}(e)_{j}d_{E}e^{j}\right) +\ \left( \sum_{j=\overline{1,m}%
}Q_{E}(e)_{j}d_{E}e^{j}\right) \wedge \left( \sum_{j=\overline{1,m}%
}Q_{E}(e)_{j}d_{E}e^{j}\right) =0,  \label{B14}
\end{equation}%
meaning that the connection $d_{E}^{\ast }:\Gamma (E^{\ast })\rightarrow
T^{\ast }(E^{\ast })\otimes \Gamma (E^{\ast })$ on $\Gamma (E)$ is flat.

\section{The loop diffeomorphisms group $\ \widetilde{Diff}(\mathbb{T}^{n}),$
the Courant type algebroid $(\mathcal{A}(\tilde{G})^{\ast },[[\cdot ,\cdot
]],\tilde{r})$ and the related integrable Hamiltonian flows}

Let us consider the product $\ \widetilde{Diff}_{+}(\mathbb{T}^{n})\times \
\widetilde{Diff}_{-}(\mathbb{T}^{n}),$ $n\in \mathbb{Z}_{+},$ where $%
\widetilde{Diff}_{\pm }(\mathbb{T}^{n})$ are subgroups of the loop
diffeomorphisms group $\ $ $\widetilde{Diff}(\mathbb{T}^{n}):=$ $\{\mathbb{C}%
\supset \mathbb{S}^{1}\rightarrow Diff(\mathbb{T}^{n})\}$ of the torus $%
\mathbb{T}^{n},$ holomorphically extended, respectively, on the interior $%
\mathbb{D}_{+}^{1}\subset \mathbb{C}$ and on the exterior $\mathbb{D}%
_{-}^{1}\subset \mathbb{C}$ \ regions of the unit centrally located disk $%
\mathbb{D}^{1}\subset \mathbb{C}^{1}$ and such that for any $\tilde{g}%
(\lambda )\in \widetilde{Diff}_{-}(\mathbb{T}^{n}),$ $\lambda \in \ \mathbb{D%
}_{-}^{1},$ $\tilde{g}(\infty )=1\in Diff(\mathbb{T}^{n}).$ The
corresponding Lie subalgebra $\widetilde{diff}(\mathbb{T}^{n})\simeq
\widetilde{diff}_{+}(\mathbb{T}^{n})\oplus \widetilde{diff}_{-}(\mathbb{T}%
^{n})\ $is a direct sum of the subalagebras $\widetilde{diff}_{\pm }(\mathbb{%
T}^{n})\simeq \widetilde{Vect}_{\pm }(\mathbb{T}^{n})\ $of the loop
subgroupss $\ \widetilde{Diff}_{\pm }(\mathbb{T}^{n})$ \ of vector fields on
$\mathbb{S}^{1}\times \mathbb{T}^{n},$ \ extended holomorphically,
respectively, on regions $\mathbb{D}_{\pm }^{1}\subset \mathbb{C}^{1},$
where for any $\tilde{a}(\lambda )\in \widetilde{diff}_{-}(\mathbb{T}^{n})$
\ the value $\tilde{a}(\infty )=0.$

Consider now the related affine Courant type algebroid $\mathcal{A}^{\ast }(%
\widetilde{Diff}(\mathbb{T}^{n}))\simeq (\widetilde{diff}(\mathbb{T}%
^{n})^{\ast }\ltimes $ $\widetilde{diff}(\mathbb{T}^{n}),\{\cdot ,\cdot \},%
\tilde{r}),$ \ where we put, by definition, the affine Lie algebra $\mathcal{%
A}(\widetilde{Diff}(\mathbb{T}^{n})):=\widetilde{diff}(\mathbb{T}%
^{n})\ltimes \widetilde{diff}(\mathbb{T}^{n})^{\ast },$ \ where the
invariant anchor morhism $\tilde{r}:\ \widetilde{diff}_{\pm }(\mathbb{T}%
^{n})^{\ast }\rightarrow \widetilde{diff}_{\pm }(\mathbb{T}^{n}),\ \tilde{r}%
:=\tilde{k}$ $\oplus \tilde{\eta},$ is such that the linear mappings
\begin{equation}
\tilde{P}_{\pm }:=I\pm \tilde{D}^{-1}/2:\mathcal{A}(\widetilde{Diff}(\mathbb{%
T}^{n}))\rightarrow \mathcal{A}_{\pm }(\widetilde{Diff}(\mathbb{T}^{n})),
\label{D12}
\end{equation}%
generated by the related derivation $\tilde{D}:=\ \tilde{k}\circ \tilde{\eta}%
^{-1}:\mathcal{A}(\widetilde{Diff}(\mathbb{T}^{n}))\rightarrow \mathcal{A}(%
\widetilde{Diff}(\mathbb{T}^{n})),$ are the corresponding projectors on the
subalagebras $\mathcal{A}_{\pm }(\widetilde{Diff}(\mathbb{T}^{n}))\subset
\mathcal{A}(\widetilde{Diff}(\mathbb{T}^{n})).$ Then the following theorem,
justifying the relationship between the Lie-Poisson bracket on the space $%
\mathcal{A}^{\ast }(\widetilde{Diff}(\mathbb{T}^{n}))$ and the corresponding
Courant structure, holds.

\begin{theorem}
The adjoint space $\mathcal{A}^{\ast }(\widetilde{Diff}(\mathbb{T}^{n}))$ is
Poissonian with respect to the following deformed Lie-Poisson structure
\begin{equation}
\{(\tilde{l},\tilde{p})(\tilde{X}),(\tilde{l},\tilde{p})(\tilde{Y})\}_{D}=(%
\tilde{l},\tilde{p})(\tilde{X},\tilde{Y}]_{D})  \label{A13}
\end{equation}%
at $(\tilde{l},\tilde{p})\in \ \mathcal{A}^{\ast }(\widetilde{Diff}(\mathbb{T%
}^{n}))$ \ for any $\tilde{X},\tilde{Y}$ $\in \mathcal{A}(\widetilde{Diff}(%
\mathbb{T}^{n})).$
\end{theorem}

Let $\ I(\mathcal{A}^{\ast }(\widetilde{Diff}(\mathbb{T}^{n})))\ :=$\ $\
\{\gamma \in D(\mathcal{A}^{\ast }(\widetilde{Diff}(\mathbb{T}^{n}))):ad_{%
\func{grad}\gamma (\tilde{l})}^{\ast }(\tilde{l},\tilde{p})=0,(\tilde{l},%
\tilde{p})\in \mathcal{A}^{\ast }(\widetilde{Diff}(\mathbb{T}^{n}))\}\ $\
denote the set of Casimir functionals on the space $\mathcal{A}^{\ast }(%
\widetilde{Diff}(\mathbb{T}^{n})),$\ invariant with respect to the canonical
co-adjoint action of the loop diffeomorphism group$\mathcal{\ }\widetilde{%
Diff}(\mathbb{T}^{n})\ $\ on $\mathcal{A}^{\ast }(\widetilde{Diff}(\mathbb{T}%
^{n})).$\ Then any countable sequence of independent functionals $\gamma
_{j}\in I(\mathcal{A}^{\ast }(\widetilde{Diff}(\mathbb{T}^{n}))),j\in
\mathbb{N},$\ are commuting to each other with respect to the deformed
Lie-Poisson structure $\{\cdot ,\cdot \}_{D}$\ \ and generate an iunfinite
hierarchy of completely integrable \ \ \cite{AbMa,BlPrSa,FaTa} \ Hamiltonian
systems
\begin{equation}
\partial (\tilde{l},\tilde{p})/\partial t_{j}=-ad_{\tilde{D}^{-1}\func{grad}%
\gamma _{j}(\tilde{l},\tilde{p})}^{\ast }(\tilde{l},\tilde{p})  \label{A14}
\end{equation}%
on the space $\mathcal{A}^{\ast }(\widetilde{Diff}(\mathbb{T}^{n}))$ with
respect to the evolution parameters $t_{j}\in \mathbb{R},j\in \mathbb{N}.$
As naturally follows from the works \ \ \cite%
{HePrBlPr,HePrBaPr-1,HePrBaPr-2,Pryk}, the constructed \ above integrable \
\ \ Hamiltonian systems \ (\ref{A14}) suitably generalize the so called \ \
\cite{Pleb} \ heavenly type differential systems, describing diverse
geometric structures of conformal type on finite dimensional Rieamnnian
manifolds.

\bigskip

\bigskip

\section{Conclusions}

There were analyzed Poisson structures related with the Courant type
algebroid, including \ those related with cotangent bundles on Lie group
manifolds. Subject to the Courant type algebroids and related R-structures \
on them there were constructed Lie-Poisson brackets invariant with respect
to the coadjoint action of the loop diffeomorphisms group and analyzed their
relationships to those generated by a special tensor mapping. \ The
corresponding integrable Hamiltonian flows, generalizing the so called
heavenly type differential systems, describing diverse geometric structures
of conformal type on finite dimensional Rieamnnian manifolds were also
described.

\bigskip

\bigskip

\section{Acknowledgements}

The authors are cordially indebted to Alexander Balinsky, Orest D.
Artemovych, Oksana Hentosh and Yarema Prykarpats for useful comments and
remarks, especially for elucidating references, which were instrumental when
preparing a manuscript. They are  also appreciated to Alina Dobrogowska for
fruitful and instructive discussions during the XL Workshop on Geometric
Methods in Physics held on 2-8 July 2023 in Bia\l owie\.{z}a, Poland. The
acknowledgments in due course belong to the Department of Computer Science
and Telecommunication of the Cracov University of Technology for a local
research grant.


\begin{thebibliography}{99}
\bibitem{AbMa} R. Abraham and J. Marsden, Foundations of Mechanics, Second
Edition, Benjamin Cummings, NY

\bibitem{Arno} V.I. Arnold. Mathematical Methods of Classical Mechanics;
Springer: Berlin/Heidelberg, Germany, 1989.

\bibitem{BlPrSa} D. Blackmore, A.K. Prykarpatsky and V.H. Samoylenko,
Nonlinear dynamical systems of mathematical physics, World Scientific
Publisher, NJ, USA, 2011

\bibitem{Blas-1} M. B\l aszak. Classical R-matrices on Poisson algebras and
related dispersionless systems. Phys. Lett. A 297(3-4) (2002), 191--195

\bibitem{Blas-2} M. B\l aszak. Multi-Hamiltonian theory of dynamical
systems. Berlin: Springer, 1998. doi: 13/978-3642637803.

\bibitem{Cour} T.J. Courant. Dirac manifolds. Transactions of the American
Mathematical Society, 319(2) (1990), 631-661

\bibitem{Mack} K.C. H. Mackenzie. Lie algebroids and Lie pseudoalgebras.
Bull. London Math. Soc. 27 (1995) 97-147

\bibitem{JiLe} V.M. Jimenez., M. de Leon. The evolution equation: an
application of grupoids to material evolution. J. Geom.\ Mech. 14(2) (2022),
331-348

\bibitem{HePrBlPr} O.Ye. Hentosh, Y.A. Prykarpatsky, Denis Blackmore, and
A.K. Prykarpatski. Dispersionless Multi-Dimensional Integrable Systems and
Related Conformal Structure Generating Equations of Mathematical Physics. \
SIGMA, 15 (2019), 079(20)

\bibitem{HePrBaPr-1} O. E. Hentosh, Ya. A. Prykarpatskyy, A. A. Balinsky,
and A. K. Prykarpatski. Geometric structures on the orbits of loop
diffeomorphism groups and related heavenly type Hamiltonian systems.: I.
Ukrainian Mathematical Journal, 74(8) (2023), 1175-1208

\bibitem{HePrBaPr-2} O. E. Hentosh, Ya. A. Prykarpatskyy, A. A. Balinsky,
and A. K. Prykarpatski. Geometric structures on the orbits of loop
diffeomorphism groups and related heavenly type Hamiltonian systems.: I.
Ukrainian Mathematical Journal, 74(8) (2023), 1175-1208

\bibitem{Pryk} A.K. Prykarpatski, Quantum Current Algebra in Action:
Linearization, Integrability of Classical and Factorization of Quantum
Nonlinear Dynamical Systems. Universe 8(2022), 288;
https://doi.org/10.3390/universe8050288

\bibitem{Godb} C. Godbillon. Geom\`{e}trie differentielle et m\'{e}canique
analytique. (French) Hermann, Paris, 1969.

\bibitem{FaTa} L.D. Faddeev, L.A. Takhtajan, Hamiltonian Methods in the
Theory of Solitons. \ Springer Berlin, Heidelberg, 2007.

\bibitem{Seme} M.A. Semenov-Tian-Shansky. What is a classical R-matrix?
Func. Anal. Appl. 1983, 17, 259--272.

\bibitem{Pleb} J. Pleba\'{n}ski, Some solutions of complex Einstein
equations. J. Math. Phys. 16 (1975), 2395--2402.

\bibitem{GaVaVo} D. Garcia-Beltran, J.A. Vallejo and Y. Vorobjev. On Lie
Algebroids and Poisson Algebras, SIGMA 8 (2012), 006(14).

\bibitem{LeMaMa} M. de Leon, J.C. Marrero, E. Marttinez. Lagrangian
submanifolds and dynamics on Lie algebroids. arXiv:math/0407528v1 [math.DG]
30 Jul 2004.

\bibitem{KoMiSl} I. Kolar, P.W. Michor, J. Slovak. Natural operations in
differential geometry, Springer-Verlag, Berlin, 1993.
\end{thebibliography}
\end{document}